\documentclass[12pt,leqno,twoside]{amsart}
\usepackage{amssymb}
\usepackage{latexsym}
\usepackage{verbatim}
\usepackage{epsfig}

\newtheorem{theorem}{Theorem}[section]
\newtheorem{corollary}[theorem]{Corollary}

\newtheorem{proposition}[theorem]{Proposition}
\theoremstyle{definition}

\theoremstyle{remark}
\newtheorem{remark}[theorem]{\sc Remark}
\theoremstyle{remark}
\newtheorem{facts}[theorem]{\sc Facts}
\theoremstyle{remark}

\theoremstyle{remark}

\theoremstyle{remark}

\theoremstyle{remark}

\renewcommand{\Box}{\square}    
\newcommand{\cal}{\mathcal}

\renewcommand{\int}{{\rm{int}}}
\newcommand{\Sing}{{\rm{Sing\hspace{1pt}}}}

\newcommand{\id}{{\rm{id}}}
\newcommand{\mult}{{\rm{mult}}}
\renewcommand{\ker}{\mathop{{\rm{ker}}}\nolimits}
\newcommand{\coker}{\mathop{{\rm{coker}}}\nolimits}

\renewcommand{\char}{{\rm{char}}}

\newcommand{\h}{{\rm{ht}}}
\newcommand{\lk}{{\rm{lk}}}
\newcommand{\e}{\varepsilon}
\newcommand{\fin}{\hspace*{\fill}$\Box$}
\newcommand{\cl}{{\rm{closure}}}

\newcommand{\cW}{{\cal W}}

\newcommand{\bC}{{\mathbb C}}

\newcommand{\bZ}{{\mathbb Z}}

\begin{document}
 
\title[Non-isolated singularities]{The vanishing neighbourhood of
non-isolated singularities}

\author{\sc Mihai Tib\u ar}

\address{Math\' ematiques, UMR 8524 CNRS,
Universit\'e des Sciences et Technologies de Lille, 
\  59655 Villeneuve d'Ascq, France.}

\email{tibar@math.univ-lille1.fr}




\begin{abstract}
We study the vanishing neighbourhood of 
 non-isolated singularities of functions on singular spaces
by associating a general linear function.
 We use the carrousel monodromy in order to show how to get a 
better control over the attaching of 
thimbles. For one dimensional singularities, we prove obstructions 
to integer (co)homology groups and 
to the eigenspaces of the monodromy via monodromies
of nearby sections. Our standpoint allows to find,
 in certain cases,  the structure of 
the Milnor fibre up to the homotopy type. 
\end{abstract}

\maketitle
\setcounter{section}{0}

\section{Introduction}\label{intro}

In the landscape of  singularities of holomorphic functions, 
the non-isolated singularities
play a particular role. Their study has been initiated in the '70 and '80 by
Y. Yomdin, R. Randell and L\^e D.T. As a natural first extension of 
isolated singularities, the case of 1-dimensional singularities got special 
attention. The first results about the homotopy type of the Milnor fibre 
were proved by D. Siersma \cite{Si-line} and his studies opened the way to a 
series of other results for 1, 2 or higher dimensional singularities,
 by R. Pellikaan,
T. de Jong, A. Zaharia, A. N\'emethi, G. Jiang, J. Fernandez and others. 
   Meanwhile, the progress in stratified Morse theory  by Goresky and
   MacPherson allowed one to treat singular holomorphic functions on singular
 spaces. The viewpoint due to  L\^e D.T. which consists in  associating
 to the function $f$ 
a general linear function $l$ and studying the couple $(l,f)$ gave rise to 
new insight in the topic of non-isolated singularities. 
  In particular, the study of the ``box'' neighbourhood lead to
new results more recently, which exploit further the properties of the 
monodromies appearing in the  fibration defined by $(l,f)$ (cf D. Siersma,
 J. Steenbrink, M. Saito,
D. Massey and others). There are some
other streams of research, which we shall not mention here. 

 Let $f:(X, 0)\to \bC$ be a
holomorphic function defined on the germ $(X, 0)$ of a singular space of pure 
dimension $n+1$, embedded into $(\bC^m,0)$ for some $m$.
 We work under the following technical but natural condition: ``the rectified 
homotopical depth of  $X$ is maximal'', which includes the cases $X$ is smooth
or a complete intersection.
 Our aim is to find how the Milnor fibre $F$ of $f$ is build from ingredients 
associated
 to restrictions of $f$ to lower dimensional slices of the space, which have
therefore lower
 dimensional singularities.
 
 It is well known that the Milnor fibre $F$ of $f$ is $(n-k)$-connected, where 
$k$ denotes the dimension of the stratified singular locus of $f$. This
follows by a Lefschetz-type argument, from
comparing $F$ to the Milnor fibre $F'$ of the restriction $f_{|l=0}$, where
$l$ is general enough.
  We go beyond this comparison and relate $F$ to $F'$ via a slice $\{l=\eta
  \}$ near to the origin (Theorem \ref{p:1}): $H_{n}(F)$ and $H_{n-1}(F)$
are the kernel and respectively the cokernel of a certain morphism $\tilde 
H_{n}(F\cup F'_D) \to H_n(F'_D, F')$, see Figure \ref{f:1} for a rapid location.
 This comparison allows us to
  exploit the $l$-monodromy and deduce bounds for the betti numbers $b_{n-1}(F)$
  and $b_{n}(F)$. The pair $(F'_D, F')$ has the advantage that it can further 
be localised, by excision, in a tubular neighbourhood of the slice 
 $\Sing f \cap F'_D$ of 
the singular locus.

 We then specialise to the case $\dim \Sing f =1$. Here we show how the 
$l$-monodromy is related to the ``vertical monodromy'' of the transversal
singularities of a slice near to the origin. We get more explicit
 bounds for $b_{n-1}(F)$ and for the eigenspaces and the maximal Jordan blocks
of the monodromy, as well as divisibility results for the 
 characteristic polynomial of
 the monodromy $h$, in terms of transversal singularities.
  These extend to a more general setting
some of Siersma's results proved in \cite{Si}. Our proof 
is based on the controlled attaching of cells, whereas Siersma's proof uses
 special variation maps.

We finally show how our standpoint allows to recover and throws a new light on
Siersma's bouquet result for ``line singularities'' (Corollary \ref{c:3,5}).  
 Further developments
 are in progress as a joint work with Dirk Siersma \cite{ST}.


\section{A geometric viewpoint}\label{geo}

 Let $(X, 0)$ be a germ of a singular space of pure dimension $n+1$ and let
 $f:(X, 0)\to \bC$ be a holomorphic function. 
We shall assume that our space $(X,0)$ satisfies one of the following 
 two conditions (where the former
implies the latter).
We refer to \cite{HL} for the notions of {\em rectified homotopical depth} $rhd(X)$
and {\em rectified homological depth} $rHd(X)$.

\smallskip 
 
\noindent {\bf (*)}  $rhd(X)\ge \dim_0 X$,

respectively

\noindent {\bf (**)}  $rHd(X)\ge \dim_0 X$.

\smallskip 

 These conditions are both true in case $X$ is a complete intersection
 (arbitrarily singular). 
 It was proved by Hamm and L\^e \cite{HL} that condition (*), resp. (**), 
implies that the complex link of $(X, p)$, for any
 $p$, 
is homotopically equivalent to a bouquet of $n$-spheres, respectively
 has the reduced homology concentrated in dimension $n$. These also imply that  
any function with a stratified Morse singularity at some point of $X$ has
 Milnor fiber homotopy equivalent to a bouquet of $n$-spheres, resp. with
 reduced cohomology concentrated in dimension $n$.  Condition $rHd(X)\ge n+1$ is further equivalent to the fact 
that  the constant sheaf on $X$ is perverse. 
 
We denote by $\Sing \phi$ the singular locus of some holomorphic map $\phi : (X,0)
\to (\bC^p,0)$, defined as the closure of the union of the singular 
loci $\Sing \phi_{|W_i}$ of the restriction of $\phi$ to
the strata  $W_i$ of some Whitney stratification $\cW$  of $X$ which we fix
throughout the paper.

Let us chose a generic linear function  $l:(X, 0)\to \bC$ and let $B_\e$
denote a Milnor ball for $f$, that is the intersection of a small 
enough ball at the origin of the ambient space with a suitable representative
of the germ $(X,0)$. L\^e D.T. showed (see e.g. \cite{Le}) that one can use
the map $(l,f) : B \to \bC^2$ and  a neighbourhood
"box" $B := B_\e \cap l^{-1}(D_{\eta'}) \cap f^{-1}(D_{\gamma'})$
in order to describe the local Milnor fibration of $f$ and its relation to the 
Milnor fibration of the slice $f_{|l=0}$.  We have already used this
 approach in \cite{Ti, Ti-iomdin} in order to get control over
 the attaching of cells via the monodromy, in case of a function $f$ with
 isolated singularity. Let us set the notations and recall several 
facts. 

 If $l$ is general enough, then the polar locus 
$\Gamma := \cl \{ \Sing(l,f)\setminus \{f=0\}\}$ is a curve on $(X,0)$, or it
is empty. Each point of the intersections $B\cap \Gamma \cap l^{-1}(\eta)$ and of $B\cap\Gamma \cap f^{-1}(\gamma)$ is a stratified Morse singularity of the restrictions $f_{l^{-1}(\eta)}$ and $l_{f^{-1}(\gamma)}$ respectively.
Let us denote by $\Delta := (l,f) (\Gamma)$ the Cerf diagram
in the target space $\bC^2$. 
Take  $\eta \ll \eta'$ and $\gamma \ll \gamma'$ such that
 the intersection $(D_\eta \times \{ \gamma\}) \cap \Delta$ with the Cerf
 diagram $\Delta$ is contained in the interior of $D_\eta \times \{ \gamma\}$.
We use the following {\bf notations and remarks}:
\begin{enumerate}
\item[(1)]  $F:= B\cap f^{-1}(\gamma)\cap l^{-1}(D_{\eta})$ is the Milnor fiber of
 $f$. Instead of the point $\gamma$ one can take any point on $\partial D_{\gamma}$.

\item[(2)]  $F' := B\cap (l,f)^{-1}(\eta, \gamma)$
  is homeomorphic to
 the Milnor fiber of the restriction $f_{|\{ l=0\}}$ (but of course has a
  different monodromy over the circle $\{\eta \} \times \partial \bar D_\eta$). 
 
\begin{figure}[hbtp]
\begin{center}
\epsfxsize=10cm
\leavevmode
\epsffile{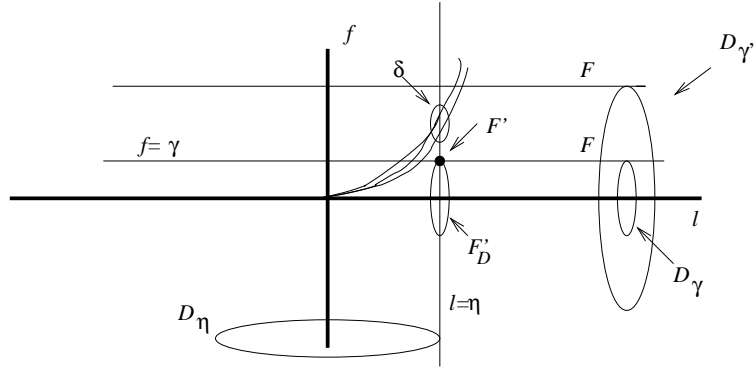}
\end{center}
\caption{{\em 
 Cerf diagram}}
\label{f:1}   
\end{figure}

\item[(3)] $F_D := B \cap f^{-1}(D_\gamma)$ is a Milnor tube, hence contractible.

\item[(4)] $F'_D:= F_D \cap l^{-1}(\eta)$ retracts to the fibre $B\cap f^{-1}(0) \cap  
l^{-1}(\eta)$,  which is the complex link of the hypersurface  $f^{-1}(0)$ at
the origin.
\end{enumerate}

We collect below a bunch of results from the literature, which will play a key
role in the following.
\begin{facts}\label{r:key}
\begin{enumerate}
\item Under condition (*), 
 the complex link of $X$ at 0 is homotopy equivalent to a
 wedge of spheres $\vee S^n$. The complex link of the hypersurface  $f^{-1}(0)$ at
the origin, denoted $\lk(f^{-1}(0),0)$, is homotopy 
equivalent to a wedge of spheres $\vee S^{n-1}$. This follows from the fact
that $f^{-1}(0)$ inherits the property (*), see \cite[Th. 3.2.1]{HL}.
Consequently $F'_D$ is homotopy 
equivalent to a wedge of spheres $\vee S^{n-1}$, since $F'_D 
\stackrel{\h}{\simeq} \lk(f^{-1}(0),0)$ .

\item The complex link of $X$ at 0, denoted $\lk(X,0)$,  is obtained by
  attaching to the complex link of $f^{-1}(0)$ at the origin, denoted
  $\lk(f^{-1}(0),0)$, a number of cells, which are all of dimension $n$
   if condition (*) is assumed. Each cell
  corresponds to a point of intersection of the polar curve $\Gamma$ with the
  slice $l^{-1}(\eta)$, see \cite{Le-cycles, Le, Ti}. 
In Figure \ref{f:1}, this can be visualised as
  attaching to $F'_D$ the cells corresponding to 
the intersection  $\Gamma\cap l^{-1}(\eta)$.

\item There is a topological disk $\{ \eta\} \times \delta$ containing  all 
intersection points $\Delta\cap \{l=\eta\}$ and is disjoint from the disk $\{
\eta\} \times D_\gamma$.
 This disk can be made sliding along the Cerf diagram to a zone $W \subset 
D_\eta \times \{ \gamma\}$.
 This procedure, due to L\^e D.T., is known as "rabattement dans le diagramme
 de Cerf" and was used in our bouquet structure theorem \cite{Ti} for the
 Milnor fibre of function germs $f$ with isolated singularity.  So,
 revisiting the point (a) above, after attaching to $F'_D$ the cells from the
 zone $(l,f)^{-1}(W)$, which we have identified to $(l,f)^{-1}(\{ \eta\}
 \times \delta)$, we get the complex link of $X$, hence a bouquet of $n$-spheres.
 If we continue to attach the other $n$-cells coming from the intersection
 points of 
 $\Gamma$ with $F$ which are outside $(l,f)^{-1}(W)$, then, as result, 
  we can only get more $n$-spheres in the bouquet. 
 \end{enumerate}
\end{facts}

Let $h$ denote the monodromy on the (co)homology of the fibre of $f$. 
This is induced by a geometric
monodromy which acts on $F$ and on $F'$, and which we shall denote by the same
symbol $h$. It acts as the identity 
on $F_D$ and on $F'_D$.

There is also the action of the geometric
 $l$-monodromy: on $F_D$, $F$ and $F'_D$ it is isotopic to the 
identity, but it may be non-trivial on $F'$.

%

\begin{theorem}\label{p:1}
Let $(X,0)$ satisfy the condition (**) and let $f: (X,0) \to \bC$ 
be any holomorphic function germ.
Then the following sequences are exact: 
\[ 0\to \tilde H_{n}(F) \to \tilde H_{n}(F\cup F'_D) \to H_n(F'_D, F')
  \to \tilde H_{n-1} (F)\to 0\]
and 
\[ 0\to \tilde H^{n-1}(F) \to H^n(F'_D, F') \to \tilde H^{n}(F\cup F'_D)\to 
\tilde H^{n} (F)\to 0\]
and there are natural isomorphisms for all $j\ge 1$:
\[ \tilde H_{n-j-1}(F) \simeq H_{n-j}(F'_D, F'), \mbox{ resp.\ }  \tilde H^{n-j-1}(F) 
\simeq H^{n-j}(F'_D, F'). \]
 
The monodromies $h$ and $l$ act on all these morphisms.

\end{theorem}
 
\begin{proof} 
Let us first consider the union $F\cup F'_D$, where $F\cap F'_D = F'$. Note
that both the monodromies $h$ and $l$ are in general not isotopic to the
identity on $F\cup F'_D$. From  Facts \ref{r:key}(b)(c)
it follows that $\tilde H_{*-1}(F\cup F'_D) = \tilde H_{*-1}(\vee  S^n)$,
and one can also deduce the exact number of spheres in the bouquet. 

Since  $F_D$ is contractible, we get $H_*(F_D, F\cup F'_D)= \tilde
H_{*-1}(F\cup F'_D)$. The same is true in cohomology.

Consider next the the pair $(F\cup F'_D, F)$. By excision we get the isomorphism:
\[ H_{*}( F\cup F'_D, F)\stackrel{\simeq}{\leftarrow} 
\tilde H_{*} (F'_D, F'),\]
and its similar counterpart in cohomology. The geometric monodromies 
 $h$ and $l$ act on all these sequences.

Let then consider the long exact sequence of the triple $(F_D, F\cup F'_D,
F)$. This splits into short sequences since
 the homology $H_*(F_D, F\cup F'_D)$ 
is concentrated in dimension $n+1$. In order to complete the proof of our
complete claim we 
just have to remark the isomorphism $H_*(F_D, F)= H_{*-1}(F)$. The proof in
 cohomology parallels the one in homology.
\end{proof}
We derive the following bound for the betti number $b_n(F)$.
\begin{corollary}\label{c:lambda}
$b_n(F) \le \lambda^0 + b_n(\lk(X,0))$, where $\lambda^0 =
  \int_0(\Gamma, f^{-1}(0)) - \int_0(\Gamma, l^{-1}(0))$.
\end{corollary}
\begin{proof}  
This follows from the first map in the above theorem and the
computation $\dim \tilde H_{n}(F\cup F'_D)= \lambda^0 + b_n(\lk(X,0))$ explained in  Facts \ref{r:key}(b)(c).
Also remember that the (co)homology of $\lk(X,0)$ is concentrated in dimension $n$, since
 condition (**) is assumed.
 Note that $\lambda^0$ depends only on 
$(X,0)$ and $f$ but not on the choice of generic $l$.
\end{proof}
 If $X$ is nonsingular, then $\lk(X,0)$ is acyclic and we recover 
the known well-known inequality $b_n(F) \le \lambda^0$.

\bigskip

We further investigate the pair $(F'_D, F')$. Let then $T$ denote a
 tubular neighbourhood of $\Sing f \cap F'_D$ within $F'_D$.
Therefore $T$ retracts to the complex link of $\Sing f$ at the origin.\footnote{In case $\Sing f$ is a complete intersection this is homotopy equivalent to a wedge of spheres of dimension  $\dim \Sing f -1$.}
Let $L$ denote the action on (co)homology of the $l$-monodromy.  We get:
\begin{corollary}\label{c:1}
\begin{enumerate}
\item  In Theorem \ref{p:1} one may replace  $H_{*}(F'_D,  F')$ by 
$H_{*}(T, T\cap F')$, and similarly in cohomology.
\item  $\coker (L-\id \mid H_{n}(T, T\cap F'))$ surjects onto $H_{n-1}(F)$, 
respectively $H^{n-1}(F)$ injects into $\ker (L-\id \mid H^{n}(T, T\cap F'))$. 
Consequently $b_{n-1}(F) \le \dim ker (L-\id \mid H^{n}(T, T\cap F'))$.
\end{enumerate}  
\end{corollary}
\begin{proof} (a) The substitution is due to the excision: $H_{*}(F'_D, F') 
\stackrel{\simeq}{\leftarrow}  H_{*}(T, T\cap F')$.
Point (b) follows from (a) and from the exact sequence of Theorem \ref{p:1}(a) 
on which
 the $l$-monodromy acts as follows. The action of the geometric
 $l$-monodromy on $F$ is isotopic to the 
identity. It is therefore the identity on the (co)homology of $F$ and of 
$(F_D, F)$. 
\end{proof}

\begin{remark}\label{r:2}
 Since the geometric $l$-monodromy is the identity on $F$, it induces the 
identity on $H_{n-j}(F'_D, F')$ and on $H_{n-j}(T, T\cap F')$
for all $j\ge 1$, by Theorem \ref{p:1}(b) and Corollary \ref{c:1}(a).
\end{remark}

\section{Case $\dim \Sing f =1$}\label{s:dim1}

We specialize to the case $\dim \Sing f = 1$ and point out several consequences of
 the preceding results. We still assume in this section that $(X, 0)$
 satisfies the condition (**). 
In case of $\dim \Sing f = 1$, the tubular neighbourhood $T$ consists of  
small Milnor balls $B_i$ at the finitely many points $\Sing f \cap F'_D$.
 Let $F_i$ denote the Milnor fiber of such an isolated singularity of the restriction of $f$
 to the transversal slice $F'_D$. Then the
geometric monodromy $h$ on $F$ restricts to the Milnor monodromy $h_i$ of
$F_i$, for each $F_i$. Let $L$ denote, as before, the action on (co)homology of the $l$-monodromy. The l-monodromy acts on the points $\Sing f \cap F'_D$ by certain permutations; each point comes back to itself after applying a number of times the l-monodromy. We denote by $\nu_i$ the action on the (co)homology of $F_i$ of the come-back monodromy\footnote{i.e. by 
definition the so-called {\em vertical monodromy}, cf \cite{St,Si}.}. One notes that
$F_i$ as well as $\nu_i$ depend only on the component $\Sigma_j$ of $\Sing f$.
 We then get:

\begin{proposition}\label{p:2}
If condition (**) holds and if $\dim \Sing f = 1$ then:
\[  \coker (L-\id\mid H_{n}(T, T\cap F')) \simeq \oplus_j\coker (\nu_j -\id \mid   
H_{n-1}(F_j)) \]
 where the sum is taken over the components $\Sigma_j$ of $\Sing f$.
\end{proposition}
\begin{proof}
 We observe that $H_{*}(T, T\cap F')$ is isomorphic, by excision, to  
$\oplus_i H_{*}(B_i,F_i) = \oplus_i \tilde H_{*-1}(F_i)$, where the sum is
taken over {\em all} the singular points which are in 
the linear slice $F'_D$. As shown in  \cite{Si-mono} or \cite{Ti-iomdin},
there is a cyclic movement of a singular point belonging to some component $\Sigma_j$ of
$\Sing f$.  This yields a particular shape of the matrix of $L$:
one can start with a basis $(e)$ of $H_{n-1}(F_1)$ and then $(L^k(e))$ is a
basis for $H_{n-1}(F_{k+1})$, for all $k\in \overline{1,s_j-1}$, where $s_j$
is the number of points of intersection of a hyperplane slice $l=\eta$ with
the component $\Sigma_j$ of $\Sing f$.
This gives the following direct sum splitting:
 \[\coker (L-\id \mid \oplus_i H_{n-1}(F_i)) \simeq \oplus_j \coker (\nu_j -\id 
\mid H_{n-1}(F_j)), \]
 where in the second sum 
we take {\em one point for each component $\Sigma_j$ of $\Sing f$}. 
\end{proof}
 \begin{corollary}\label{c:2}
 \begin{enumerate}
\item $b_{n-1}(F) \le \dim \oplus_{\Sigma_j} \coker (\nu_j -\id \mid  H_{n-1}(F_j))$.
  
\item 
 $\char_{h |H_{n-1}(F)}$ divides $\prod_{\Sigma_j} \char_{h_j|\coker (\nu_j -\id\mid
 H_{n-1}(F_j))}$.\\
 In particular $\char_{h |H_{n-1}(F)}$ divides the product
 $\prod_{\Sigma_j} \char_{h_j |H_{n-1}(F_j)}$. 
\end{enumerate} 
\end{corollary}
 In case of non-singular $X$, this result was proved by Dirk Siersma \cite{Si}
 with a different proof.
 A weaker version, yet for singular $(X,0)$ satisfying
(**), was proved in \cite{Ti-iomdin}.
\begin{proof}  From Corollary \ref{c:1}(b) and Proposition \ref{p:2} we deduce the surjection: 
\[ \oplus_j \coker (\nu_j -\id\mid H_{n-1}(F_j))\twoheadrightarrow  
 H_{n-1}(F), \]
from which (a) follows immediately.\\
(b). The monodromy $h$ acts on the preceding surjection. On the left hand side this amounts to the 
action of $\oplus_j h_j$. One remarks that the monodromies 
  $h_j$ and $\nu_j$ commute. 
Then apply $\char_h$ to this
and get the claimed divisibility. 
\end{proof}
Let us remark that the analogous results are true in cohomology by
standard reasons; 
one just replaces ``surjection'' by ``injection'' and ``$\coker$'' by
 ``$\ker$''.

\bigskip
Let $b_\lambda(V, \mu)$ denote the dimension of the eigenspace
corresponding
to the eigenvalue $\lambda$ of the linear operator $\mu$ acting on the vector space
$V$. Let  $J_\lambda(V, \mu)$ denote the maximum of the sizes of the 
Jordan blocks. With these notations we have:

 \begin{corollary}\label{c:3}
In cohomology, let $K_j :=  \ker (\nu_j -\id \mid  H^{n-1}(F_j))$. Then:
 \begin{enumerate}
\item $H^{n-1}(F) \subset \oplus_{\Sigma_j} K_j$.
\item $b_\lambda(H^{n-1}(F), h) \le \sum_{\Sigma_j} b_\lambda(K_j, h_j)$.
\item $J_\lambda(H^{n-1}(F), h) \le \sum_{\Sigma_j} J_\lambda(K_j, h_j)$.

\end{enumerate} 
\end{corollary}
\begin{proof} (a) is clear from the preceding remark. \\
(b) and (c). The monodromy $h$ acts on the inclusion morphism (a),
so $h - \lambda \id$ acts too.  On the right hand side term, this amounts to 
the action of $h_i$, respectively of $h_i - \lambda \id$, on $H^{n-1}(F_i)$,
independently for each $i$. We know that $h_i$ commutes with $\nu_i$.
\end{proof}
  The study may be  pursued  in case 
of higher dimensional singular locus $\Sing
f$; this is work in progress jointly with Dirk Siersma \cite{ST}.



\section{Controlled attaching}\label{s:attach}
This is a study of 
the cell-attaching described at the end
of Facts \ref{r:key}(b) 
in the case $\dim \Sing f =1$ and $(X,0) = (\bC^{n+1},0)$.
We refer to the notations and results in \S \ref{geo}; let us 
recall that condition (*) is fulfilled in this case.

We shall use, and therefore need to recall from \cite{Ti} the construction of
 relative thimbles associated to the pair $(F,F')$ which are adapted to the 
carrousel monodromy. The disk $D_\eta \times \{ \gamma \}$ will be called
 {\em carrousel disk}.

 Our study on the attaching of thimbles starts with the following 5 steps. 
We refer to Figure \ref{f:1}. The first 3 steps apply to a singular space
$(X,0)$ under the condition (*).

\smallskip
\noindent
(1). Our assumptions imply that  $F'$ is a bouquet of $n-1$ spheres. The restriction
 of $f$ to the slice $l=\eta$ has 
  only isolated singularities: the intersections with $\Sing f$
  and the intersections with $\Gamma$ (which are of Morse
 type, by the genericity of $l$).  Therefore this is a deformation of the 
singularity $f_{|l=0}$ which is a partial Morsification. 

\smallskip
\noindent
(2). The result of attaching to $F'$ the 
thimbles corresponding to the singularities in the zone $(l,f)^{-1}(\{ \eta\} 
\times \delta)$ is, up to homotopy type, a bouquet of spheres of dimension
$n-1$. The number of spheres is equal to the sum of the Milnor numbers of the
singularities of the slice $\{ f=0\}\cap \{ l= \eta\}$.
 Now, remember from Facts \ref{r:key}(b) that we may identify the zone
 $(l,f)^{-1}(\{ \eta\} \times \delta)$ to $(l,f)^{-1}(W)$, where $W$ is a
 certain 
 open subset of the
 carrousel disk $D_{\eta}\times \{\gamma\}$ defined in \cite[\S 2]{Ti}.
 
The Milnor fiber 
$F := (l,f)^{-1}(D_\eta\times \{ \gamma \})$ is obtained from $F':= 
(l,f)^{-1}(\eta, \gamma)$ by
attaching a number of $n$-cells, each cell corresponding to one of the
 intersection points $\Gamma \cap F$. The total number of cells is
  equal to the intersection multiplicity
$\mult_0(\Gamma, \{ f=0\})$. 
We  first attach to $F'$ the thimbles from the zone
$(l,f)^{-1}(W)$. We have seen before what is the result of this attaching.

\smallskip
\noindent
(3). The {\em carrousel model of the monodromy}, introduced by L\^e D.T.\cite{Le},
together with the refined description of the carrousel monodromy from
\cite{Ti}, see 
\cite[Fig. 1, pag. 233]{Ti}, enables one to describe 
how the further attaching occurs.

If $\Delta$ is not empty, 
then there exist thimbles
  outside the zone $(l,f)^{-1}(W)$; 
this is due to the fact
 that 
 all the components of $\Delta$ are tangent to the horizontal axis $\{ f=0\}$.

 Let us attach one ``next'' thimble, which is out of the zone $(l,f)^{-1}(W)$).
 By the main construction in \cite{Ti}, this thimble is the image by the 
{\em carrousel monodromy} of a 
thimble from the zone $(l,f)^{-1}(W)$. Call the latter thimble $t^0$ and the
 former $t^1$.
 Say $t^0$ attaches to $F'$ over the cycle $a$. By the construction in
 \cite{Ti}, the thimble $t^1$
 will attach to $F'$ exactly over the cycle $h_1(a)$, where $h_1$ denotes the geometric
 $f$-monodromy in the slice $F'_D$ around the singular points $\Sing f \cap
 \{l=\eta\}$.
  All the above explanation is already contained in \cite{Ti}, where the case $\dim
 \Sing f= 0$
is treated; in that case the monodromy $h_1$ is geometrically trivial.

\smallskip
\noindent 
(4). What is $h_1(a)$ more precisely in the case $(X,0) = (\bC^{n+1},0)$ and
 $\dim \Sing f= 1$? We consider the cycles up to homotopy 
equivalence
and denote by $[\cdot ]$ the homology or homotopy equivalence classes. In particular
 the additive notation for homotopy classes means that our homotopy groups are
 restricted to  dimension $\ge 2$ and $\ge 3$ for relative homotopy groups
 (i.e. thimbles) respectively.

 Let us remark, firstly,  that the $f$-monodromy in the slice $F'_D$ around the
 singular points $\Sing f \cap \{l=\eta\}$ splits into the direct sum of 
the monodromies around each of these points, since they are all in the same
fibre of $f$. Secondly, that each such
singularity may be more complex than a Morse singularity and so 
the monodromy around each of such
singular points is the Coxeter element of some Morsification of the
corresponding singularity. 
Then, by applying the Picard-Lefschetz rules to the monodromy $h_1$ we get: 
\begin{equation} \label{eq:piclef}
h_1([a]) = [a] + \sum_r k_r [b_r],
\end{equation}
where $k_r$ is an integer and $b_r$ denotes one of the cycles vanishing 
at some point of $\Sing f \cap
\{l=\eta\}$ (and where the sum is taken over these cycles).

\smallskip
\noindent
(5). After all this discussion, we come back to the attaching of $t^0$ and
$t^1$. Let $a$ be as chosen before.
Notice now that the attaching of $t^0$ over $a$ just kills $a$, in other words
$a$ is contractible in the space $F'\cup (l,f)^{-1}(W)$.
 Next, the attaching of $t^1$ to $F'\cup (l,f)^{-1}(W)$ is done over a cycle
 homotopy equivalent to $[a] + \sum_r k_r [b_r]$. Since $a$ is 
 contractible in the space $F'\cup (l,f)^{-1}(W)$, this new attaching is 
really over the cycle $\sum_r k_r [b_r]$. We have shown before 
 that this cycle is not zero, and it is clearly
 independent on the cycles of type $a$.
We get the following conclusion:

\begin{proposition}\label{p:attach}
Let $(X,0) = (\bC^{n+1},0)$ and
 $\dim \Sing f= 1$.
Assume that the linear function $l$ is sufficiently general. Then $F$ is obtained from $F'\cup (l,f)^{-1}(W)$ by
 attaching thimbles defined by the carrousel monodromy, with attaching maps
of the type $\sum_r k_r [b_r]$, where $b_r$'s are the cycles of $F'$ vanishing
at the singular points $\Sing f \cap \{l=\eta\}$.\fin
\end{proposition} 
An exceptional case is when no thimbles are attached to $F'\cup (l,f)^{-1}(W)$.  
This is equivalent to $\Delta$ being not tangent to the axis $f=0$,
which can only happen if  $\Delta = \emptyset$ for generic $l$, in which 
case one has $F\stackrel{\h}{\simeq} F'$.  This is in turn equivalent 
 to the fact that the singular locus
 $\Sing f$ is a line (by the non-splitting principle of L\^e D.T.) and that the
   Milnor number of the transversal
singularity is constant along the line.  The 
equivalences can be deduced from the attaching results discussed in Facts
\ref{r:key}.  

\smallskip

Let us further remark that there exists a geometric cycle of type $a$ (i.e. a cycle over
which  attaches one of the carrousel thimbles from the zone $(l,f)^{-1}(W)$) 
which has
 non-zero intersection with  at least one of the $b_r$'s. This is due to the
 fact that the cycles of types $a$ and $b$ are together a basis of cycles (which 
 one may arrange to be a geometric basis) 
 of the isolated singularity at the origin of $f_{| l=0}$. Their  intersection graph is
connected, cf \cite{La}. It is then
 an easy exercise* (which we may safely leave to the reader) to show that there
 exists at least 
one $[b_r]$ in the sum (\ref{eq:piclef}) such that its
 coefficient $k_r$ is non-zero. 
 
 The attaching map is precisely
the boundary map $\partial : H_{n}(F,F') \to H_{n-1}(F')$ in homology. In homotopy, the
additive notation has a meaning only if $n\ge 3$ and then the attaching map is the corresponding boundary morphism $\partial : \pi_{n}(F,F';.) \to \pi_{n-1}(F';.)$.
Now remark that the monodromy $h$ acts on the exact sequence of the pair $(F, F')$, 
and its action on $H_{n-1}(F')$ or $\pi_{n-1}(F')$ is precisely $h_1$.
Therefore $h_1(\partial (\alpha))= \partial(h(\alpha))$ for any linear combination $\alpha$ of carrousel thimbles. This proves the following, where coefficients are in some field:

 \begin{corollary}\label{c:4}
 Under the hypothesis of Proposition \ref{p:attach},  assume in addition that
 the singular locus
 $\Sing f$ is not a line with constant transversal Milnor number. Then
the betti number $b_{n-1}(F)$ is strictly less than the sum $\sum_i b_{n-1}(F_i)$ of Milnor numbers 
of the singularities in the slice $\{ l= \eta\} \cap \{ f=0\}$. \fin
\end{corollary}
This corollary has been announced independently by L\^e D.T.
and D.B. Massey \cite{LM}. 

\noindent
It is necessary to compare this result to Siersma's one in \cite{Si} which
coincides to Corollary \ref{c:2}(a) in case of a non-singular 
space germ $(X, 0)$. When $\Sing f$ is a
union of lines, then Corollary \ref{c:4} provides a bound which might be
  better by at most one.  Nevertheless in all other cases Siersma's result
 gives a better bound, or at least the same. 

\begin{corollary}\label{c:3,5}
Let $f$ define a line singularity with A$_1$-transversal generic singularity
 type and
 different from the transversal type at the origin. Let $[b]$ be the cycle
 vanishing at the A$_1$-transversal generic singularity.
Then:
\begin{enumerate}
\item In homology, there exists a linear combination of carrousel thimbles
 which attaches to $F'\cup (l,f)^{-1}(W)$ 
with attaching map $[b]$. 
\item The Milnor fibre $F$ of $f$ is homotopy equivalent to 
a bouquet $\vee S^{n}$. 
 \end{enumerate}
\end{corollary}
\begin{proof}
Let $[a_i]$ denote the cycles of $F'$ which are killed by carrousel thimbles from 
the zone $(l,f)^{-1}(W)$. 
We claim that there exists a linear combination of carrousel thimbles
$[t_i]$ such that the attaching map to $F'$ is:
 \[  \sum_i s_i h_1(t_i^0) \mapsto  [b] + \sum_i s_i [a_i],\]
where $s_i\in \bZ$. 

The claim follows from the following  fact: if  $p_i = \langle a_i, b\rangle$
denotes the intersection pairing
 then $\gcd \{ |p_i |\}_i =  1$.
Indeed, if $v:= \gcd \{p_i\}_i > 1$ then all $p_i$'s are zero modulo $v$,
which contradicts the fact that the Dynkin diagram of $f_{|l=0}$ relative 
to a distinguished basis is
connected in the homology with coefficients any $\bZ$-module (see
e.g. \cite[p. 77]{AGV}). 

Since $\gcd \{ |p_i |\}_i =  1$,  there exist integer coefficients $s_i$ such that
$\sum_i s_i p_i = 1$. Therefore $\langle \sum_i s_i [a_i], [b]\rangle =  1$.
By Proposition \ref{p:attach} the linear combination of thimbles 
$\sum_i s_i h_1(t_i^0)$ attaches to $F'\cup (l,f)^{-1}(W)$ with attaching map 
$\langle \sum_i s_i [a_i], [b]\rangle [b]$, which is
homotopy equivalent to $[b]$. So
this attaching will kill the cycle $[b]$.

\noindent
(b). Part (a) implies that the reduced homology of $F$ is concentrated in
dimension $n$. For $n\ge 3$ the above attaching result holds in homotopy too.
 We shall treat the cases $n=1, 2$ separately.

 For any $n\ge 2$ and in particular for $n=2$, we may apply a result by
 L\^e-Saito \cite{LS}: this says that in our situation the fundamental group 
of $F$ is abelian, hence trivial (since $H_1(F) =0$). Then $F$ has the homotopy groups of 
a bouquet $\vee S^n$, via the Hurewicz map, and therefore, by the Whitehead 
theorem for CW-complexes,  it is homotopy
equivalent to $\vee S^n$.
As for the case $n =1$, $F$ has the homotopy type of a connected
$1$-dimensional CW-complex, hence $F\stackrel{\h}{\simeq}\vee S^1$.
\end{proof}
Corollary \ref{c:3,5}(b) recovers, with a different proof,
 Siersma's bouquet result \cite{Si-line}.


\bigskip

\noindent
(*)  {\em Solution of the exercise.} 
Take as the order of the cycles indexed by $r$'s 
  the order of composing them in the Coxeter element.
 Then to the first index $r$ such
that $\langle a, b_r\rangle  \not= 0$ will correspond the coefficient $k_r = \pm 
\langle a, b_r\rangle$. The following $k_r$'s are, in general, linear
 combinations of intersection numbers.

\end{document}